\documentclass[12pt]{amsart}
\title[Airy-heat, Hermite and higher order Hermite]{Airy-heat functions, Hermite and higher order Hermite generating functions}
\author{Gerardo Hern\'andez-del-Valle}
\begin{document}
\maketitle
\begin{abstract}
In this note we discuss the relationship between the generating functions of some Hermite polynomials $H$,
\begin{eqnarray*}
\sum\limits_{j=0}^\infty H_{j\cdot n}(u)\frac{z^n}{n!},
\end{eqnarray*}
 generalized Airy-Heat equations
 \begin{eqnarray*}
 \frac{1}{2\pi}\int_{-\infty}^{+\infty}\exp\left\{a(i\lambda)^n-\frac{1}{2}\lambda^2t+i\lambda x\right\}d\lambda,
 \end{eqnarray*}
 higher order PDE's
 \begin{eqnarray*}
 \frac{\partial u}{\partial t}(t,x)=a\frac{\partial^nu}{\partial x^n}(t,x)+\frac{1}{2}s\frac{\partial^2u}{\partial x^2}(t,x)
 \end{eqnarray*}
  and generating functions of higher order Hermite polynomials $H^{(n)}$:
  \begin{eqnarray*}
  \sum\limits_{j=0}^\infty H^{(n)}_j(v)\frac{x^j}{j!}.
  \end{eqnarray*}
  In particular, we show that under some conditions, these problems are equivalent.
  \end{abstract}
  \section{Introduction}
  This note is motivated by the observation, made in Hern\'andez-del-Valle (2010), that all the known closed form densities of hitting times of Brownain motion may be expressed in terms of heat polynomials, which may be alternatively formulated as Hermite polynomials. 
  
  Except in the case of the constant and linear boundaries, it seems that solutions of the hitting time problem will involve infinite series of Hermite polynomials. In general, these series expansions will not be convergent, as has been pointed out in Widder and Rosenbloom (1959). 
  
  Hence, in the present work we derive an alternative representation, which relates to generalized Airy functions, higher order Hermite polynomials, and higher order PDE's.

  \section[Triple lacunary and Airy-heat equations]{Triple Lacunary generating function of Hermite polynomials and Airy-Heat equation}
  In Gessel and Jayawant (2005) the authors derive---using umbral and combinatorial arguments---a generating function for Hermite polynomials of order 3, namely:
  \begin{eqnarray}\label{hermite}
  \sum\limits_{n=0}^\infty H_{3n}(u)\frac{z^n}{n!}=\frac{e^{(w-u)(3u-w)/6}}{\sqrt{1-6wz}}\phantom{1}_2F_0\left(\frac{1}{6},\frac{5}{6};-:\frac{54z^2}{(1-6wz)^3}\right)
  \end{eqnarray}
 in which $w=(1-\sqrt{1-12uz})/6z=u\cdot C(3uz)$, where $C(x)=(1-\sqrt{1-4x})/2x$ is the Catalan number generating function, and $\!\!\phantom{1}_2F_0$ is the so-called `mysterious' hypergeometric function, which in turn is related to the Bessel function $B$ of order $1/3$ through
 \begin{eqnarray*}
 \phantom{1}_2F_0\left(\frac{1}{6},\frac{5}{6};-;\xi\right)=\sqrt{-\frac{1}{\xi^2}}e^{-\frac{1}{2\xi^2}}B_{1/3}\left(-\frac{1}{2\xi^2}\right),
 \end{eqnarray*}
 and alternatively, $B_{1/3}$ is (under certain conditions) equivalent to the so-called Airy function Ai, defined as:
 \begin{eqnarray}\label{airy}
 \hbox{Ai}(x):=\frac{1}{2\pi}\int_{-\infty}^{+\infty}\exp\left\{i\lambda +i\lambda x\right\}d\lambda,
 \end{eqnarray}
 where $i:=\sqrt{-1}$.
\subsection{Alternative derivation of (\ref{hermite})} A solution to the {\it forward\/} heat equation
\begin{eqnarray*}
\frac{\partial w}{\partial t}(t,x)=\frac{1}{2}\frac{\partial^2w}{\partial x^2}(t,x)
\end{eqnarray*}
is given in terms of the following {\it inverse Fourier transform}:
\begin{eqnarray}\label{ai1}
w(t,x)=\frac{1}{2\pi}\int_{-\infty}^{+\infty}e^{i\frac{\lambda^3}{3}-\frac{1}{2}\lambda^2 t+i\lambda x}d\lambda.
\end{eqnarray}
After completing the cube in the exponential function, equation (\ref{ai1}) may be written as:
\begin{eqnarray*}
w(t,x)=\exp\left\{\frac{t^3}{12}+\frac{tx}{2}\right\}\hbox{Ai}\left(x+\frac{t^2}{4}\right)
\end{eqnarray*}
where Ai is defined in (\ref{airy}).

Finally, equations (\ref{hermite}) and (\ref{ai1}) are `symbolically' equivalent, by interchanging the McLaurin expansion of $e^{i\lambda^3/3}$ in (\ref{ai1})
\begin{eqnarray}
\label{aa} w(t,x)&=&\frac{1}{2\pi}\int_{-\infty}^{+\infty}e^{i\frac{\lambda^3}{3}-\frac{1}{2}\lambda^2 t+i\lambda x}d\lambda\\
\nonumber&=&\frac{1}{2\pi}\int_{-\infty}^{+\infty}e^{-\frac{1}{2}\lambda^2 t+i\lambda x}\sum\limits_{j=0}^n\frac{1}{j!}\left(i\frac{\lambda^3}{3}\right)^jd\lambda\\
\label{new}&=&\sum\limits_{j=0}^n\frac{(-3)^{-j}}{j!}\left[\frac{1}{2\pi}\int_{-\infty}^{+\infty}e^{-\frac{1}{2}\lambda^2 t+i\lambda x}(-i\lambda)^{3j}d\lambda\right]
\end{eqnarray}
and observing that the term within the squared brackets correspond to the so-called {\it derived\/} heat polynomials of order $n$---introduced by Widder and Rosenbloom (1959):
\begin{eqnarray}
\nonumber\omega_n(t,x)&:=&\frac{1}{2\pi}\int_{-\infty}^{+\infty}e^{i\lambda x-\frac{1}{2}\lambda^2 t}(-i\lambda)^nd\lambda\\
\label{triple}&=&t^{-n/2}k(t,x)H_n\left(\frac{x}{\sqrt{2t}}\right).
\end{eqnarray}
It follows that (\ref{new}), together with (\ref{triple}) is equivalent to
\begin{eqnarray}\label{triple}
w(t,x)=\sum\limits_{j=0}H_{3n}(u)\frac{z^n}{n!}.
\end{eqnarray}

As a corollary, we conclude that generating functions of Hermite polynomials of the form
\begin{eqnarray*}
\sum\limits_{j=0}^\infty H_{j\cdot m}(u)\frac{z^j}{j!}
\end{eqnarray*}
are equivalent to Airy-heat functions, defined as:
\begin{eqnarray*}
\frac{1}{2\pi}\int_{-\infty}^{+\infty}\exp\left\{a\lambda^m-\frac{1}{2}\lambda^2 t+i\lambda x\right\}d\lambda,
\end{eqnarray*}
with the appropriate identification of variables $u$ and $z$.

There is though an important issue, the interchange of summation and integration in (\ref{new}) is only symbolic, since $e^{i\lambda y^3/3}$---although entire---is not of order $(2,\sigma)$. In other words, approximating (\ref{aa}) with (\ref{new}) is not allowed.
\section{Correspondence between expansions of Hermite and Higher order Hermite polynomials}
In this section we describe a sort of  `duality' between expansions of the form:
\begin{eqnarray*}
\sum\limits_{j=0}^\infty H_{j\cdot n}(u)\frac{z^n}{n!}
\end{eqnarray*}
 and generating functions of higher order Hermite polynomials $H^{(n)}$:
 \begin{eqnarray*}
 \sum\limits_{j=0}^\infty H^{(n)}_j(v)\frac{x^j}{j!}.
 \end{eqnarray*}
defined as:
\begin{eqnarray*}
H^{(n)}_j(t,x)&=&e^{t\partial^n_x}x^j\\
&=&j!\sum\limits_{k=0}^{[j/n]}\frac{x^{j-nk}t^k}{(j-nk)!k!}.
\end{eqnarray*}

But first, let us embed solutions of the heat equation into higher order PDE's:
\subsection{Higher order PDEs and the heat equation}
Let us start by defining the following bivariate PDE
\begin{eqnarray}\label{high}
\frac{\partial u}{\partial t}(t,x)=a\frac{\partial^n u}{\partial x^n}(t,x)+\frac{1}{2}s\frac{\partial^2u}{\partial x^2}(t,x)
\end{eqnarray}
where $a$ and $s$ are arbitrary coefficients. It follows immediately that it has particular solutions such as:
\begin{eqnarray*}
u(t,x)&=&\exp\left\{\left(a\lambda^n+\frac{1}{2}s\lambda^2\right)t+\lambda x\right\}\\
&=&\exp\left\{\left(a(i\lambda)^n-\frac{1}{2}s\lambda^2\right)t+i\lambda x\right\}\\
&=&\exp\left\{\left(a(i\lambda)^n-\frac{1}{2}s\lambda^2\right)t\right\}\sin(\lambda x)\\
&=&\exp\left\{\left(a(i\lambda)^n-\frac{1}{2}s\lambda^2\right)t\right\}\cos(\lambda x)
\end{eqnarray*}
or
\begin{eqnarray}\label{sep}
u(t,x)=e^{\lambda t}v(x)
\end{eqnarray}
where $v$ is a solution to the following ODE
\begin{eqnarray*}
\lambda v(x)=av^{(n)}(x)+\frac{1}{2}sv^{(2)}(x).
\end{eqnarray*}
Which follows after substitution of (\ref{sep}) in (\ref{high}). 

Alternatively, applying the Fourier transform with respect to the space variable $x$ to (\ref{high}):
\begin{eqnarray*}
\tilde{u}_t&=&a(i\lambda)^n\tilde{u}-\frac{1}{2}s\lambda^2\tilde{u}\\
&=&\left[(ai\lambda)^n-\frac{1}{2}s\lambda^2\right]\tilde{u},
\end{eqnarray*}
and solving for $t$, leads to
\begin{eqnarray*}
\tilde{u}(t,\lambda)=\tilde{u}(0,\lambda)\exp\left\{a(i\lambda)^nt-\frac{1}{2}s\lambda^2t\right\}.
\end{eqnarray*}
Hence, in terms of the {\it inverse\/} Fourier transform equation (\ref{high}) admits solutions of the following form
\begin{eqnarray}\label{heathigh}
u(t,x)=\frac{1}{2\pi}\int_{-\infty}^{+\infty}\tilde{u}(0,\lambda)\exp\left\{a(i\lambda)^nt-\frac{1}{2}s\lambda^2t+i\lambda x\right\}d\lambda.
\end{eqnarray}

For example, in solution (\ref{heathigh}) let $a=-1/3$ and $u(0,\lambda)=1$ 
\begin{eqnarray*}
u(t,x)=\frac{1}{2\pi}\int_{-\infty}^{+\infty}\exp\left\{\left(i\frac{\lambda^3}{3}-\frac{1}{2}s\lambda^2\right)t+i\lambda x\right\}
\end{eqnarray*}
which after completing the cube is equivalent to
\begin{eqnarray}\label{ai}
u(t,x)=\exp\left\{\frac{s^3}{12}t+\frac{sx}{2}\right\}t^{-1/3}\hbox{Ai}\left(t^{-1/3}\left\{x+\frac{s^2t}{4}\right\}\right).
\end{eqnarray}
In particular, when $t=1$, then (\ref{ai}) is a solution to the heat equation in variables $(s,x)$ and alternatively when $s=0$, it leads to an Airy kernel, introduced by Widder (1975).

Note that, either by a proper time change or simply by setting $t=1$ in (\ref{heathigh})
\begin{eqnarray*}
u(1,x)=\frac{1}{2\pi}\int_{-\infty}^{+\infty}\tilde{u}(0,\lambda)\exp\left\{a(i\lambda)^n-\frac{1}{2}s\lambda^2+i\lambda x\right\}d\lambda
\end{eqnarray*}
leads to solutions of the heat equation in terms of variables $(s,x)$:
$$\frac{\partial v}{\partial s}(s,x)=\frac{1}{2}\frac{\partial^2 v}{\partial x^2}(s,x).$$
On the other hand if we set $s=0$ and $\tilde{u}(0,\lambda)=1$ in (\ref{heathigh}) 
\begin{eqnarray*}
u(t,x)=\frac{1}{2\pi}\int_{-\infty}^{+\infty}\exp\left\{a(i\lambda)^nt+i\lambda x\right\}d\lambda
\end{eqnarray*}
we get solutions to higher order PDE's which may be alternatively described in terms of higher order `Airy' functions:
\begin{eqnarray*}
Ai^{(n)} (u):=\frac{1}{2\pi}\int_{-\infty}^{+\infty}\exp\left\{a(i\lambda)^{3+n}+i\lambda x\right\}d\lambda
\end{eqnarray*}
(that is, when $n=0$ we obtain the classic Airy function), i.e.
\begin{eqnarray}\label{airydensity}
u(t,x)=\frac{1}{t^{1/(n+3)}}Ai^{(n)}\left[\frac{x}{t^{1/(n+3)}}\right].
\end{eqnarray}

In general, with $\tilde{u}(0,\lambda)=1$,  it follows that (\ref{heathigh}) is a convolution of (\ref{airydensity}) and
the heat kernel:
\begin{eqnarray}\label{hk}
v(s,x)=\frac{1}{\sqrt{2\pi ts}}\exp\left\{-\frac{x^2}{2st}\right\}.
\end{eqnarray}
We may now describe the duality between generating functions of Hermite and higher order Hermite polynomials.
\subsection{Duality between generating functions of Hermite and higher order Hermite polynomials}
As described in the previous subsection a solution to
\begin{eqnarray*}
\frac{\partial u}{\partial t}(t,x)=a\frac{\partial^n u}{\partial x^n}(t,x)+\frac{1}{2}s\frac{\partial^2u}{\partial x^2}(t,x)
\end{eqnarray*}
is
given by
\begin{eqnarray*}
u(t,x)=\frac{1}{2\pi}\int_{-\infty}^{+\infty}\exp\left\{a(i\lambda)^nt-\frac{1}{2}s\lambda^2t+i\lambda x\right\}d\lambda.
\end{eqnarray*}
which can be alternatively seen as the convolution 
\begin{eqnarray}\label{con}
u(t,x):=\int_{-\infty}^{+\infty}v(\tau,y)w(t,x-y)dy
\end{eqnarray}
between the fundamental solutions of the following pair of PDE's:
\begin{eqnarray}
\label{u}\frac{\partial w}{\partial t}(t,x)&=&\frac{\partial^n w}{\partial x^n}(t,x)\\
\nonumber\frac{\partial v}{\partial \tau}(\tau,x)&=&\frac{\partial^2v}{\partial x^2}(\tau,x)
\end{eqnarray}
[see equations (\ref{airydensity}) and (\ref{hk}) respectively] where we have set $a=1$ and $\tau=t\cdot s$.

On the other hand recall the definition of the higher-order Hermite polynomials:
\begin{eqnarray}\label{genh}
H^{(n)}_j(t,x)&=&e^{t\partial^n_x}x^j
\end{eqnarray}
which after differentiation with respect to $t$ on both sides of (\ref{genh}) we have that
\begin{eqnarray*}
\partial_t H^{(n)}_j(t,x)&=&\partial^n_x\left[e^{t\partial^n_x}x^j\right]\\
&=&\partial^n_xH^{(n)}_j(t,x)
\end{eqnarray*}
is a solution of (\ref{u}) and at time $t=0$
$$H_j(n)(0,x)=x^j,$$
i.e. it solves a Cauchy problem. After substiution of  $v$ equation (\ref{hk}) in (\ref{con}) and subsequently taking its MacLaurin expansion, it follows:
\begin{eqnarray*}
u(t,x)&:=&\int_{-\infty}^{+\infty}\frac{1}{\sqrt{2\pi\tau}}\exp\left\{-\frac{y^2}{2\tau}\right\}w(t,x-y)dy\\
&=&\frac{1}{\sqrt{2\pi\tau}}\sum\limits_{j=0}^\infty\left(-\frac{1}{2\tau}\right)^j\frac{1}{j!}\int_{-\infty}^{+\infty}y^{2j}w(t,x-y)dy\\
&=&\frac{1}{\sqrt{2\pi\tau}}\sum\limits_{j=0}^\infty\left(-\frac{1}{2\tau}\right)^j\frac{1}{j!}H^{(n)}_{2j}(t,x).
\end{eqnarray*}
So for instance, a symbolically equivalent triple lacunary generating function described in (\ref{triple}) based upon higher order Hermite polynomials is given by:
\begin{eqnarray*}
\frac{1}{\sqrt{2\pi t}}\sum\limits_{j=0}^\infty\left(-\frac{1}{2t}\right)^j\frac{1}{j!}H^{(3)}_{2j}(1,x).
\end{eqnarray*}

\end{document}